\documentclass[11pt]{amsart}
\usepackage{amssymb,amsbsy,upref,epsf,MnSymbol}
\usepackage{amsmath,amssymb,amsthm}
\usepackage{graphicx,mathtools,mathrsfs,enumitem}
\usepackage[margin=1in]{geometry}
\usepackage{fancyhdr}
\usepackage{hyperref}
\newtheorem{theorem}{Theorem}
\newtheorem{proposition}[theorem]{Proposition}
\newtheorem{lemma}[theorem]{Lemma}
\newcommand{\R}{\mathbb{R}}

\newcommand{\norm}[1]{\left\lVert#1\right\rVert}
\newcommand{\paren}[1]{\left( #1 \right)}
\newcommand{\bracket}[1]{\left[ #1 \right]}
\newcommand{\abs}[1]{\left\lvert #1 \right\rvert}

\DeclareMathOperator{\sgn}{sign}
\newcommand{\del}{\partial}

\newcommand{\ddt}{\frac{d}{dt}}

\newcommand{\rest}{\big|}

\newcommand{\n}{^{-1}}
\newcommand{\eps}{\varepsilon}
\newcommand{\indic}[1]{\chi_{\{#1\}}}
\newcommand{\spo}{s_-}
\newcommand{\sne}{s_+}

\title[Contraction of Viscous Shocks]{$L^2$-type Contraction of Viscous Shocks for Large Family of Scalar Conservation Laws}

\author[Stokols]{Logan F. Stokols} 
\address[L. F. Stokols]{\newline Department of Mathematics, \newline The University of Texas at Austin, Austin, TX 78712, USA}
\email{lstokols@math.utexas.edu}

\date{\today}


\thanks{\textbf{Acknowledgment.} This work was partially supported by the NSF Grant DMS 1614918.}

\begin{document}

\begin{abstract}
In this paper we study small shocks of 1D scalar viscous conservation laws with uniformly convex flux and nonlinear dissipation.  We show that such shocks are $L^2$ stable independent of the strength of the dissipation, even with large perturbations.  The proof uses the relative entropy method with a spatially-inhomogeneous psuedo-norm.  
\end{abstract}

\maketitle \centerline{\date}
This paper will consider 1D scalar dissipative conservation equations of the form
\begin{equation} \label{eq:main} 
\del_t u + \del_x\left[ Q(u) \right] = \nu \del_{xx} \eta'(u), 
\end{equation}
where $Q$ and $\eta$ are uniformly convex functions, meaning that for some constant $\Lambda \geq 1$,
\begin{equation} \label{convexity condition}
\frac{1}{\Lambda} \leq \eta''(x), Q''(x) \leq \Lambda
\end{equation}
holds for all $x \in \R$.  This bound on $\eta''$ is natural because $\eta''$ measures the coercivity of the dissipation in divergence form.

Equations of this form admit a class of traveling wave solutions known as shocks.  Shocks are monotone decreasing and exponentially constant at $\pm \infty$.  Given any two values $\spo > \sne$, there exists a shock $s:\R \to (\sne, \spo)$ such that
\begin{align*}
\lim_{s \to \infty} s(x) &= \sne, \\
\lim_{s \to -\infty} s(x) &= \spo,
\end{align*}
and $s(x-t \sigma)$ is a solution to \eqref{eq:main} with constant
\[ \sigma := \frac{Q(\spo) - Q(\sne)}{\spo - \sne}. \]
This formula for $\sigma$ is known as the Rankine-Hugoniot condition.  Viscous shocks are a generalization of inviscid shocks, which are piece-wise constant with a single jump discontinuity.  Inviscid shocks are recovered in the limit as $\nu \to 0$.  

We will show in this paper that sufficiently small shock solutions are $L^2$-stable.  Since even small perturbations in $L^2$ can significantly affect the travelling speed of a shock, we will show stability only up to a Lipschitz shift which depends on the perturbation.  This limitation is not present in the $L^1$ theory (see Kruzkhov \cite{Kr.entropy}), but is well known in the theory of $L^2$ shock stability (see Leger \cite{Le}).  

We will prove the following:
\begin{theorem}[Main Theorem] \label{thm:main simple}
Let $\Lambda \geq 1$ be a constant, and let $\eta, Q: \R \to \R$ be satisfy \eqref{convexity condition} on the interval $(-R,R)$ for $R \in (0,+\infty]$, and let $\eta'''$, $Q'''$ continuous at 0.  Then there exists a constant $\eps_0$ such that the following holds:

Let $\nu > 0$ be any constant.  Let $s: \R \times [0,\infty) \to [\sne,\spo]$ be a shock solution to \eqref{eq:main} with $|\sne - \spo| = 2\eps \leq 2\eps_0$, and let $u:\R\times[0,\infty)\to\R$ be a solution to \eqref{eq:main} such that $\norm{u(\cdot,0) - s(\cdot,0)}_{L^2(\R)} < \infty$.  If $R < \infty$, assume $\norm{u(\cdot,0)}_{L^\infty(\R)} < R$.  

Then there exists a Lipschitz function $\gamma:[0,\infty) \to \R$ such that for any $t \in [0,\infty)$ we have
\[ \int |u(x,t) - s(x-\gamma(t),t)|^2 \,dx \leq 4\Lambda^2 \int |u(x,0) - s(x,0)|^2 \,dx. \]

The quantity $\norm{\gamma'}_{L^\infty}$ depends only on $\eps$, $\Lambda$, and $\norm{u(\cdot,0) - s(\cdot,0)}_{L^2(\R)}$.
\end{theorem}
Notice that this result is independent of the strength $\nu$ of the dissipation.  

We prove this result using the method of relative entropy, first introduced by DiPerna and Dafermos \cite{Da} to study stability of Lipschitz solutions of conservation laws.  This method has since been applied by Vasseur, Serre, Leger, and others (\cite{SeVa}, \cite{LeVa}, \cite{Le}) to show $L^2$ stability of shocks under large perturbations.  

For an entropy function $f$, we denote the relative entropy between two solutions $u_1$ and $u_2$ by
\[ f(u_1|u_2) := f(u_1) - f(u_2) - f'(u_2) [u_1-u_2]. \]
In this paper, we will use the function $\eta$, the antiderivative of the dissipative term, as our entropy function.  Our proof involves taking the time-derivative of the relative entropy of $u$ relative to the shock $s$.  Because of the assumption \eqref{convexity condition}, the integral of the relative entropy is essentially equal to the $L^2$ norm.  However, this quantity will not decrease in general, as shown by Vasseur and Kang in \cite{KaVa.burgers}.  We supplement the method by considering a weighted psuedo-norm, as in \cite{Va.entropy} and \cite{Va.recent}.  The weight function $a$ is independent of solution $u$, and is approximately constant.  

We will show the following result, from which Theorem~\ref{thm:main simple} follows as a corollary:
\begin{theorem} \label{thm:main technical}
Let $\Lambda \geq 1$ be a constant, and let $\eta, Q: \R \to \R$ satisfy \eqref{convexity condition} for all $x \in \R$, and let $\eta'''$, $Q'''$ continuous at 0.  Let $\nu = 1$.  Then there exists a constant $\eps_0$ such that the following holds:

Let $0 < \eps < \eps_0$ be a constant and let $s: \R \to [\sne,\spo]$ be a stationary shock solution to \eqref{eq:main} with $s_\pm = \mp \eps$.  Then there exists a weight function $a:\R \to [1/2,2]$ such that the following holds:

For any $u:\R\times[0,\infty)\to\R$ solving \eqref{eq:main} such that $\norm{u(\cdot,0) - s(\cdot)}_{L^2(\R)} < \infty$, there exists a Lipschitz function $\gamma:[0,\infty) \to \R$ such that for any $t \in [0,\infty)$ we have
\[ \ddt \int a(x + \gamma(t)) \eta(u(x,t)|s(x+\gamma(t))) \,dx \leq -\eps_0 \int a(x+\gamma(t)) \abs{\del_x(\eta'(u)-\eta'(s))}^2 \,dx. \]

The quantity $\norm{\gamma'}_{L^\infty}$ depends only on $\eps$, $\Lambda$, and $\norm{u(\cdot,0) - s(\cdot)}_{L^2(\R)}$, and $\norm{a-1}_\infty$ tends to 0 as $\eps \to 0$.  
\end{theorem}

The theory of $L^2$ stability of shocks is contrasted with the $L^1$ theory, as in the work of Kruzkov \cite{Kr.entropy}.  See also Ilyin and Oleinik \cite{IlOl} and Freistuhler and Serre \cite{FrSe}.  Unlike Kruzkov's result, we only need one entropy.  Though 1D scalar laws have infinitely many entropies in general, systems of conservation laws typically only have one entropy so methods which rely on multiple entropies are more difficult to generalize, though such generalizations exist, see for example Bressan, Liu, and Yang \cite{BrLiYa}.  The $L^p$ stability theory has also been studied by Adimurthi, Ghoshal, and Veerappa Gowda \cite{AdGhVG}.  $L^2$ stability has been studied outside the context of relative entropy, as by Goodman \cite{Go}, though wish stronger assumptions on the perturbation.  

Since our result is independent of the strength $\nu$ of dissipation, it is well suited to taking an inviscid limit.  

The technique used in this paper has previously been applied by Kang and Vasseur to certain 1D dissipative systems in \cite{KaVa.navier} (including 1D isotropic Navier-Stokes) and 1D scalar equations with constant dissipation in \cite{Ka} (i.e. $\eta'(u)=u$).  We are able to consider arbitrary convex dissipation by utilizing $\eta$ as an entropy.  

As in \cite{KaVa.navier}, the proof proceeds by braking up the solution $u$ into a part which is $L^\infty$ close to $s$ and an error term which may be large in $L^\infty$.  The close part is handled similarly to the existing literature, while for the error term we must make careful use of the relationship between the dissipative term and the derivative of the weight function $a$.  

The paper is structured as follows: in Section~\ref{sec:prelim} we compute the time derivative of the relative entropy.  In Section~\ref{sec:lemmas} we present a number of lemmas which will be used throughout the paper.  In Section~\ref{sec:functional} we show that our expression for the derivative of the relative entropy is non-positive under a number of special assumptions.  Finally in Section~\ref{sec:main thm} we prove Theorems~\ref{thm:main simple} and \ref{thm:main technical}.  

\section{Time derivative} \label{sec:prelim}

For any function $f$, define
\[ f(x|y) := f(x) - f(y) - f'(y)(x-y),\]
In particular, for $\eta$ our entropy the quantity $\eta(u|s)$ is called the relative entropy of $u$ relative to $s$.  

We call $\eta$ an entropy function if there exists a function $G$ such that
\[ G'(x) = Q'(x) \eta'(x).\]
In the 1D case, such a $G$ trivially exists.  

We also define
\[ F(x;y) := G(x) - G(y) - \eta'(y)\left[ Q(x) - Q(y) \right].\]

We begin by computing the time derivative of the relative entropy with arbitrary shift and arbitrary weight. 

\begin{proposition}[Time Derivative]\label{thm:time derivative}
Let $u:\R \times [0,\infty) \to \R$ and $s:\R \to \R$ be solutions to \eqref{eq:main} with $\nu=1$ and $s$ a stationary solution.  Assume $\abs{u(\cdot,t) - s(\cdot)} \in L^2(\R)$ for all $t$.  

Then for any differentiable function $\gamma:[0,\infty) \to \R$ and weight function $a \in L^\infty(\R)$, we have
\[ \ddt \int a(x+\gamma(t)) \eta\big(u(x)\big|s(x+\gamma(t)\big) \,dx = R(u) := \dot{\gamma} Y(u) + B(u) - D(u) \]
where
\begin{align*}
Y(u) &:= \int a' \eta\left(u|s\right) dx - \int a s' \eta''(s)(u-s) dx, \\
D(u) &:= \int a |\del_x (\eta'(u) - \eta'(s))|^2 dx, \\
B(u) &:= \int a' F(u;s) dx - \int a \eta''(s)s'Q(u|s) dx + \int \frac{a''}{2} (\eta'(u) - \eta'(s))^2 dx  + \int a \eta'(u|s) \del_x Q(s) dx.
\end{align*}

Here it is understood that $u$ is evaluated always at $(x,t)$ while $a$ and $s$ are evaluated at $x+\gamma(t)$.  
\end{proposition}

The expressions $Y(u)$, $B(u)$, $D(u)$ will be referenced throughout this paper with the definitions given above, and they will be abbreviated as $Y$, $B$, $D$ when the input $u$ is clear from context.  

\begin{proof}
Initially, we have
\[ \ddt \int a(x+\gamma(t)) \eta\left(u(x)|s(x+\gamma(t)\right) \,dx = \int a'\dot{\gamma} \eta\left(u|s\right) \,dx + \int a \left[ \eta'(u) - \eta'(s) \right] \del_t u \,dx \]
\[ + \int a \left[ -\eta''(s) (u-s) \right] \del_t s \,dx. \]
Since $\del_t s = \dot{\gamma} s'$, we have, with $Y$ as defined in the theorem statement, 
\begin{equation} \label{time deriv with Y}
\ddt \int a \eta(u|s) \,dx = \dot{\gamma} \bracket{Y} + \int a \left[ \eta'(u) - \eta'(s) \right] \del_t u \,dx 
\end{equation}

Note that
\[ \del_t u = \del_{xx} \eta'(u) - \del_x Q(u) \]
and that because $s$ is a shock with zero drift,
\[ \del_x Q(s) = \del_{xx} \eta'(s). \]
Therefore, writing $w = \eta'(u)-\eta'(s)$, 
\begin{equation}\begin{aligned} \label{expansion of deltu term}
\int a w \del_t u \,dx &= \int a w \bracket{ \del_{xx} \eta'(u) - \del_x Q(u) } \,dx + \int a \bracket{ -\eta''(s) (u-s) } \left[ \del_{xx} \eta'(s) - \del_x Q(s) \right] \,dx
\\ &= \int a w \del_{xx} w \,dx + \int a \eta'(u|s) \del_{xx} \eta'(s) \,dx - \int a \paren{w \del_x Q(u) - \eta''(s) (u-s) \del_x Q(s)} \,dx
\end{aligned}\end{equation}

Now, notice that
\begin{equation}\begin{aligned} \label{F term explanation}
\del_x F(u;s) + \eta''(s)Q(u|s)s' &= \bracket{ \eta'(u) Q'(u) - \eta'(s)Q'(u) } \del_x u + \bracket{ -\eta''(s)(Q(u)-Q(s)) } s' \\ & \hspace{1cm} + \eta''(s) Q(u|s) s'
\\ &= \bracket{ \eta'(u) - \eta'(s) } Q'(u) \del_x u + \eta''(s)s' \bracket{ Q(u|s) - \paren{Q(u)-Q(s)} }
\\ &= \bracket{ \eta'(u) - \eta'(s) } \del_x Q(u) - \eta''(s) (u-s) Q'(s) s'
\\ &= w \del_x Q(u) - \eta''(s) (u-s) \del_x Q(s).
\end{aligned}\end{equation}

Combining \eqref{time deriv with Y}, \eqref{expansion of deltu term}, and \eqref{F term explanation} we obtain
\[ \ddt \int\!\! a \eta(u|s) \,dx = \dot{\gamma}\bracket{Y}
+ \int\!\! a w \del_{xx} w \,dx 
+ \int\!\! a \eta'(u|s) \del_x Q(s) \,dx
- \int\!\! a \left[ \del_x F(u;s) + \eta''(s)s'Q(u|s) \right] \,dx. \]

Integrating by parts, we have
\[ \int a w \del_{xx} w \,dx = \frac{1}{2} \int a'' w^2 \,dx - \int a |\del_x w|^2 \,dx \]
\[ = \frac{1}{2} \int a'' w^2 \,dx - D \]
and
\[ \int a \del_x F(u;s) \,dx = - \int a' F(u;s) \,dx. \]

The proposition follows.  

\end{proof}


Notice that each term in $Y$ and $B$ contain either a derivative of $s$ or a derivative of $a$.  This inspires us to choose our weight function $a$ to be a linear transformation of $s$.  We can then perform a change of variables and simplify the expression even further.  The new variable $y = \eta'(s(x))$ is known as the entropic variable.  

\begin{lemma} \label{thm:rewrite time derivative}
Under the same assumptions as Proposition \ref{thm:time derivative}, if $a := 1 - \frac{\lambda}{\eps} \eta'(s)$ for some $\lambda > 0$ then, in terms of the variable $y := \eta'(s)$, we have
\begin{align*} 
Y &= \frac{\lambda}{\varepsilon} \int \eta(u|s) \,dy + \int a (u-s) \,dy, \\
D &= -\int a \eta''(s) s' |\del_y w|^2 \,dy, \\
B &= \frac{\lambda}{\eps} \int F(u;s) \,dy + \int a Q(u|s) \,dy + \frac{\lambda}{\eps}\int \frac{Q'(s)}{2\eta''(s)} w^2 \,dy - \int a \frac{Q'(s)}{\eta''(s)} \eta'(u|s) \,dy. 
\end{align*}
\end{lemma}

\begin{proof}
Notice first that $x \mapsto \eta'(s)$ is a monotone-decreasing differentiable bijection, so $u$ is a well-defined function of $y$.  The integrating factor for this new variable is 
\[ dy = -\eta''(s)s' \,dx. \]
Note the minus sign because $s'$ is negative so the direction of integration is reversed.  

The derivatives of $a$ are
\[ \del_x a = - \frac{\lambda}{\eps} \eta''(s) s' \]
and
\[ \del_{xx}a = - \frac{\lambda}{\varepsilon} \del_{xx} \eta'(s) = - \frac{\lambda}{\varepsilon} \del_x Q(s). \]
The derivative of $Q(s)$ is
\[ \del_x Q(s) = Q'(s) s' = - \frac{Q(s)}{\eta''(s)} \eta''(s) s'. \]

From here, the form of $Y$ and $B$ are trivial to compute.  

For $D$, we must simply compute
\[ \del_x w = \eta''(s) s' \del_y w. \]
\end{proof}

\section{Lemmas} \label{sec:lemmas}
This section consists of a series of lemmas which will be necessary throughout the rest of the paper.  

We begin by applying Taylor's formula to each of the quantities appearing in the expressions $Y(u)$, $B(u)$, and $D(u)$ defined in Lemma~\ref{thm:rewrite time derivative}.  These estimates, together with the bounds on the derivatives of $\eta$ and $Q$, will be the basis of all our control on the quantities $Y$, $B$, $D$.  

\begin{lemma}\label{thm:taylor}
Let $x_1$ and $x_2$ be real numbers.  Then the following estimates hold:
\begin{enumerate}[label=(\alph*)]
\item \label{taylor u-s}
There exists a point $z_0$ between $x_1$ and $x_2$ such that
\[ x_1-x_2 = \frac{1}{\eta''(z_0)}(\eta'(x_1)-\eta'(x_2)). \]
\item \label{taylor Q}
There exists a point $z_1$ between $x_1$ and $x_2$ such that
\[ Q(x_1|x_2) = \frac{Q''(z_1)}{2\eta''(z_0)^2}(\eta'(x_1)-\eta'(x_2))^2. \]
\item \label{taylor eta}
There exists a point $z_2$ between $x_1$ and $x_2$ such that
\[ \eta(x_1|x_2) = \frac{\eta''(z_2)}{2\eta''(z_0)^2} (\eta'(x_1)-\eta'(x_2))^2. \]
\item \label{taylor eta' thrd}
There exists a point $z_3$ between $x_1$ and $x_2$ such that
\[ \eta'(x_1|x_2) = \frac{\eta'''(z_3)}{2\eta''(z_0)^2} [\eta'(x_1)-\eta'(x_2)]^2. \]
\item \label{taylor eta' scnd}
There exists a point $z_4$ between $x_1$ and $x_2$ such that
\[ \eta'(x_1|x_2) = \paren{1 - \frac{\eta''(x_2)}{\eta''(z_4)}}[\eta'(x_1) - \eta'(x_2)]. \]
\item \label{taylor F}
There exists a point $z_5$ between $x_1$ and $x_2$ such that
\[ F(x_1;x_2) = \frac{1}{2}\eta''(z_5) \frac{Q'(z_5)}{\eta''(z_0)^2} (\eta'(x_1)-\eta'(x_2))^2 + \frac{1}{2}\eta'(z_5)-\eta'(x_2)] \frac{Q''(z_5)}{\eta''(z_0)^2} (\eta'(x_1)-\eta'(x_2))^2. \]
\item \label{taylor 1-y2}
If $s$ is a stationary shock solution to \eqref{eq:main} with $\nu = 1$, and $\varsigma \in (\sne,\spo)$ is a real number, then there exist points $z_6, z_7, z_8 \in (\sne,\spo)$ such that
\[ -\eta''(s) s'\rest_{s=\varsigma} = \frac{Q''(z_6)}{2\eta''(z_7) \eta''(z_8)} [\eta'(\varsigma) - \eta'(\sne)][\eta'(\spo) - \eta'(\varsigma)]. \]
\end{enumerate}
\end{lemma}

\begin{proof}
Claim \ref{taylor u-s} follows immediately from Taylor's theorem:
\[ \eta'(x_1) = \eta'(x_2) + \eta''(z_0) (x_1-x_2). \]

Applying Taylor's theorem to $Q$,
\[ Q(x_1) = Q(x_2) + Q'(x_2)(x_1-x_2) + \frac{Q''(z_1)}{2} (x_1-x_2)^2. \]
Therefore 
\[ Q(x_1|x_2) = \frac{Q''(z_1)}{2} (x_1-x_2)^2 \]
and \ref{taylor Q} follows from \ref{taylor u-s}.  

Claims \ref{taylor eta} and \ref{taylor eta' thrd} follow by the same logic as \ref{taylor Q}.  

Apply \ref{taylor u-s} to the definition of $\eta'(x_1|x_2)$ to obtain
\[ \eta'(x_1|x_2) = [\eta'(x_1)-\eta'(x_2)] - \frac{\eta'(x_2)}{\eta'(z_0)} [\eta'(x_1)-\eta'(x_2)] \]
and \ref{taylor eta' scnd} follows.  

For \ref{taylor F}, we can calculate, by Taylor's theorem,
\begin{align*} 
F(x_1; x_2) &= F(x_2;x_2) + \frac{d}{dx_1} F(x_2;x_2) (x_1-x_2) + \frac{1}{2} \frac{d^2}{dx_1^2} F(t_5;x_2) (x_1-x_2)^2 
\\ &= 0 + 0 + \frac{1}{2} \bracket{\eta''(t_5) Q'(t_5) - Q''(t_5) \bracket{\eta'(t_5) - \eta'(x_2)}} (x_1-x_2)^2. 
\end{align*}
From this and \ref{taylor u-s}, the claim \ref{taylor F} follows.  

\vskip.3cm

Since $s$ is a shock solution, $\del_x Q(s) = \del_{xx} \eta'(s)$.  Moreover $\del_x \eta'(s)\rest_{s=\spo} = 0$.  Therefore
\[ -\eta''(\varsigma) s'\rest_{s=\varsigma} = -\del_x \eta'(\varsigma)\rest_{s=\varsigma} = Q(\spo) - Q(y). \]
Now since $Q(\sne) = Q(\spo)$ by the Rankine-Hugoniot condition, there exists a point $z_6 \in (\sne,\spo)$ such that
\[ Q(\spo) - Q(\varsigma) = Q''(z_6) (\varsigma-\sne)(\spo-\varsigma). \]
Applying \ref{taylor u-s} a final time, the proof is complete.   
\end{proof}

The following lemma is Proposition 3.3 in \cite{KaVa.navier}.  It is a Poincar\'{e} type inequality.  
\begin{lemma}[Poincar\'{e}] \label{thm:poincare}
Given a constant $C_1$, there exists a constant $\delta_0 > 0$, such that for any $\delta \leq \delta_0$ the following holds:

For any $W \in L^2(0,1)$ such that $\sqrt{x(1-x)} \del_x W \in L^2(0,1)$ with $\norm{W}_2^2 \leq C_1$, the quantity
\[ \frac{-1}{\delta}\left( \int_0^1\!\! W^2 dx + 2 \!\!\int_0^1\!\! W dx \right)^2 + (1+\delta) \!\!\int_0^1\!\! W^2 dx + \frac{2}{3}\int_0^1\!\! W^3 dx + \delta \!\!\int_0^1\!\! |W|^3 \,dx - (1-\delta) \!\!\int_0^1\!\! x(1-x)|\del_x W|^2 dx\]
is non-positive.  

\end{lemma}

The following lemma is a kind of weighted Gagliardo-Nirenberg interpolation.  The quantity $D(u)$ defined in Lemma~\ref{thm:rewrite time derivative} controls the second derivative of $w$ but that control degenerates near the endpoints.  The lemma interpolates between $D$ and the $L^2$ norm to control arbitrary $L^p$ norms.  

\begin{lemma}[Gagliardo-Nirenberg] \label{thm:gagliardo}
Let $h > 0$, $p\geq 1$, $L > 0$, and $\bar{C} \leq 2 h^2 L$ be constants.  For any $w \in L^2([-L,L])$ with
\[ \int_{-L}^L w^2 \,dy \leq \bar{C}, \]
define
\[ \tilde{D} := \int_{-L}^L (y - L)(L-y) \indic{|w|>h} \abs{\del_y w}^2 \,dy. \]

Then for any $q \in (0,1)$ there exists a constant $C_q$ depending only on $q$ such that
\[ \int (w-h)_+^p \,dy \leq C_q \paren{h^{-2} \bar{C}}^q |L|^{-p/2} \tilde{D}^{p/2}. \]
\end{lemma}

\begin{proof}
By Chebyshev's inequality, $|\{|w|>h\}| \leq h^{-2} \bar{C}$ so since $h^{-2} \bar{C} \leq L$ there exists a point $y_0 \in [-L/2,L/2]$ such that $(w-h)_+(y_0) = 0$.  

For any other point $y_1$, we can calculate
\begin{align*}
|(w-h)_+(y_1)| &= |(w-h)_+(y_1) - (w-h)_+(y_2)| 
\\ &\leq \int_{y_0}^{y_1} |\del_y (w-h)_+| \,dy
\\ &\leq \paren{\int_{y_0}^{y_1} \bracket{(L+y)(L-y)}^{-1} \,dy}^{1/2} \paren{\int_{y_0}^{y_1} \bracket{(L+y)(L-y)} \abs{\del_y (w-h)_+}^2 \,dy}^{1/2}
\\ &\leq \paren{ \frac{1}{2L} \bracket{\ln(L+y) - \ln(L-y)}_{y_0}^{y_1} }^{1/2} \tilde{D}^{1/2}
\\ &= \frac{\tilde{D}^{1/2}}{(2L)^{1/2}} \bracket{\ln(L+y_1) - \ln(L-y_1) - \ln(L+y_0) + \ln(L-y_0)}^{1/2}
\end{align*}
Since $y_0 \in [-L/2,L/2]$, we can estimate $\frac{L+y_0}{L-y_0} \in [1/3,3]$ so $\ln(L+y_0) - \ln(L-y_0)$ is bounded.  The expression $\ln(L+y_1) - \ln(y_1-L)$ is similarly bounded for $|y_1|<L/2$.  For $|y_1| > L/2$, the $\ln(L-y_1)$ term will dominate for $y_1$ positive and the $\ln(y_1+L)$ term will dominate for $y_1$ negative, so for some constant $C$ we have the bound
\[ |(w-h)_+(y_1)| \leq C \paren{\frac{\tilde{D}}{L}}^{1/2} \max\paren{1,\abs{\ln(L-|y_1|)}}^{1/2}. \]

Let $\mu \leq h^{-2}\bar{C}$ be the measure of the set $\{|w|>h\}$.  Without loss of generality we assume that this region is concentrated near $\pm L$, and so 
\begin{align*}
\int (w-h)_+^p \,dy &\leq 2 \int_{L-\mu/2}^L C^p \paren{\frac{\tilde{D}}{L}}^{p/2} \abs{\ln(L-|y_1|)}^{p/2} \,dy
\\ &\leq C \paren{\frac{\tilde{D}}{L}}^{p/2} \int_0^{\mu/2} |\ln(x)|^{p/2} \,dx
\\ &\leq C_q \paren{\frac{\tilde{D}}{L}}^{p/2} \mu^q.
\end{align*}
Here we have used an estimate of the integral of $\ln(x)$ near the origin which uses the fact that $\ln(x)$ grows slower than any power of $x$.  

Since $\mu \leq h^{-2}\bar{C}$, the lemma follows.  
\end{proof}

The following final lemma shows that the quantity $Y$ bounds the $L^2$ norm.  
\begin{lemma}\label{thm:Y bounds L2}
There exists a constant $C = C(\Lambda)$ so that the following holds:

Let $\eta$ and $Q$ as in Theorem~\ref{thm:main technical} and $u$,$s$ be any functions such that $u-s \in L^2(\R)$.  Let $Y(u)$ be as in Lemma~\ref{thm:rewrite time derivative}. Then the function $w := \eta'(u)-\eta'(s)$ satisfies
\[ \int w^2 \,dy \leq C(\Lambda) \bracket{\frac{\eps}{\lambda} |Y(u)| + \frac{\eps^3}{\lambda^2}}. \]
\end{lemma}

\begin{proof}
From the definition of $Y$, we know that
\[ \frac{\lambda}{\eps} \int \eta(u|s) \,dy \leq |Y| + \int a (u-s) \,dy. \]
The right-hand side is of course non-negative since $\eta$ convex.  

Recall the notation $w = \eta'(u) - \eta'(s)$.  From Lemma~\ref{thm:taylor} \ref{taylor eta} and \ref{taylor u-s} we know that $\eta(u|s) \geq \Lambda^{-3}  w^2$ and $|u-s| \leq \Lambda |w|$.  Of course $|a| \leq 2$.  Therefore

\[ \int w^2 \,dy \leq \Lambda^3 \frac{\eps}{\lambda} |Y| + \Lambda^3 \frac{\eps}{\lambda} 2 \Lambda \int |w| \,dy. \]

By H\"{o}lder's inequality, $2 \int |w| \,dy \leq \frac{\lambda}{2 \Lambda^4 \eps} \int w^2 \,dy + \frac{2 \Lambda^4 \eps}{\lambda} \int 1 \,dy$.  Thus
\[ \int w^2 \,dy \leq \Lambda^3 \frac{\eps}{\lambda} |Y| + \frac{1}{2} \int w^2 \,dy + 2 \Lambda^8 \frac{\eps^2}{\lambda^2} \int 1 \,dy. \]

Since 
\[ \int 1 \,dy = \eta'(\spo)-\eta'(\sne) \leq 2 \Lambda \eps, \]
the lemma follows.  
\end{proof}

\section{Functional Estimates} \label{sec:functional}
In this section, we consider the quantity $-Y(u)^2 + B(u) - D(u)$ under certain assumptions on $u$.  Note that we do not need to assume $u$ is a solution of \eqref{eq:main} in this section at all, only that $u$ and $s$ are in some sense small functions.  

\begin{proposition}[Decrease for small perturbations] \label{thm:functional near}
Let $\eta$ and $Q$ satisfy \eqref{convexity condition} for all $x \in \R$ and have $\eta'''$, $Q'''$ continuous at 0.  For any positive constant $\bar{C}$, there exist constants $h_1 > 0$ and $\eps_1 > 0$, such that the following holds:

Let $s$ be a stationary shock solution to \eqref{eq:main} with $\nu = 1$ and $s_{\pm}= \mp \eps$ with $0 < \eps < \eps_1$, and let $\bar{u} \in L^\infty(\R)$ be such that $|\bar{w}| := |\eta'(\bar{u}) - \eta'(s)| \leq h$ for some $0 < h < h_1$.  Let $0 < \lambda < \eps_1$ and $a := 1-\frac{\lambda}{\eps} \eta'(s)$ such that $1/2 \leq a \leq 2$. Assume
\[ \int \bar{w}^2 \leq \bar{C} \frac{\eps^3}{\lambda^2}. \]
Then
\[ \bar{R} := \frac{-1}{2\eps^2} Y(\bar{u})^2 + B(\bar{u}) - (1-h) D(\bar{u}) \]
is non-positive.  
\end{proposition}

In the case that $\eta$ and $Q$ are quadratic polynomials, for example if $\Lambda=1$, this theorem would hold by a straightforward application of Lemma~\ref{thm:poincare}. Since $\eta$ and $Q$ have continuous second derivatives, for small inputs their second derivatives will be nearly constant and we can treat them as polynomials.  We will use Taylor's theorem, specifically in the form of Lemma~\ref{thm:taylor}, to formalize this observation.  

\begin{proof}
Let $\delta_0$ be the constant indicated by Lemma~\ref{thm:poincare} corresponding to constant $\Lambda \bar{C}$, and consider arbitrary $0 < \delta \leq \delta_0$.  

We will estimate $Y$, $B$, and $D$ using the formulae provided in Lemma~\ref{thm:rewrite time derivative}.  Notice that, since $\eta'''$ and $Q'''$ exist and are continuous at 0, $\eta''$ and $Q''$ must also be continuous at 0.  

First we analyze the term $Y$.  Define
\[ Y_1 := \frac{\lambda}{\eps} \int \eta(\bar{u}|s) \,dy. \]
By Lemma~\ref{thm:taylor} \ref{taylor eta}, there exist $t_1$, $t_2 \in [-\eps_1-h_1, \eps_1+h_1]$ so 
\[ \abs{ Y_1 - \frac{1}{2 \eta''(0)} \frac{\lambda}{\varepsilon} \int \bar{w}^2 \,dy } = \abs{\frac{\eta''(t_1)}{2\eta''(t_2)^2} - \frac{1}{2\eta''(0)}} \frac{\lambda}{\varepsilon} \int \bar{w}^2 \, dy. \]
Since $\eta''$ is continuous at 0, for $\eps_1 + h_1$ sufficiently small we can say
\[ \abs{ Y_1 - \frac{1}{2 \eta''(0)} \frac{\lambda}{\varepsilon} \int \bar{w}^2 \,dy } \leq \delta \frac{\lambda}{\eps} \int \bar{w}^2 \,dy. \] 

Define
\[ Y_2 = \int a (u-s) \,dy \]
and, by applying Lemma~\ref{thm:taylor} \ref{taylor u-s}, we can argue as above that for $\eps_1 + h_1$ sufficiently small we have
\begin{align*}
\abs{Y_2 - \eta''(0)\n \int \bar{w} \,dy } &= \int \bracket{\eta''(t_1)\n a - \eta''(0)\n} \bar{w} \,dy
\\ &= \int \eta''(t_1)\n (a-1) \bar{w} \,dy + \int \bracket{\eta''(t_1)\n - \eta''(0)\n} \bar{w} \,dy 
\\ &\leq \int \paren{\lambda \Lambda + \abs{\eta''(t_1)\n - \eta''(0)\n}} |\bar{w}| \,dy
\\ &\leq C (\lambda + \delta) \int |\bar{w}| \,dy.
\end{align*}

Since $Y = Y_1 + Y_2$, assuming without loss of generality $\eps < \delta$, we can apply the general formula $-(a+b)^2 \leq -\paren{\frac{\eps}{2\delta}} a^2 + \frac{\eps}{\delta} b^2$ for $a,b \in \R$ and $\eps/\delta \in (0,1]$ to obtain
\[ -Y^2 \leq \frac{\eps}{8 \delta} \eta''(0)^{-2} \paren{\frac{\lambda}{\eps} \int \bar{w}^2 \,dy + 2 \int \bar{w} \,dy}^2 + C\frac{\eps}{\delta} \paren{(\lambda + \delta) \int |\bar{w}| \,dy + \delta \frac{\lambda}{\eps} \int \bar{w}^2 \,dy}^2 \]
Since $\int |\bar{w}| \,dy \leq C \eps^{1/2} \sqrt{ \int \bar{w}^2 \,dy}$ and $\int \bar{w}^2 \,dy \leq \bar{C} \eps^3/\lambda^2$,
\begin{equation}\begin{aligned} 
- Y^2 &\leq \frac{\eps}{\delta} \frac{\Lambda^2}{8} \paren{\frac{\lambda}{\eps} \int \bar{w}^2 \,dy + 2 \int \bar{w} \,dy}^2 
	+ C(\lambda + \delta)^2 \frac{\eps}{\delta} \paren{\int |\bar{w}| \,dy}^2 + C\frac{\eps}{\delta} \delta^2 \frac{\lambda^2}{\eps^2} \frac{\eps^3}{\lambda^2} \int \bar{w}^2 \,dy
\\ &\leq \frac{\eps}{\delta} \frac{\Lambda^2}{8} \paren{\frac{\lambda}{\eps} \int \bar{w}^2 \,dy + 2 \int \bar{w} \,dy}^2 
	+ C\paren{\frac{\eps^2 \lambda^2}{\delta} + \eps^2\delta + \delta} \int \bar{w}^2 \,dy. 
\label{Y2(ubar) bound}
\end{aligned} \end{equation}

\vskip.3cm

Now we analyze $B$.

For the relative flux term, we estimate by Lemma~\ref{thm:taylor} \ref{taylor Q} and continuity of $\eta''$ and $Q''$
\begin{equation} \label{Q(ubar) bound} \begin{aligned}
\abs{\int a Q(\bar{u}|s) \,dy - \frac{Q''(0)}{2\eta''(0)^2} \int \bar{w}^2 \,dy} 
&\leq \int \abs{a \frac{Q''(t_1)}{2\eta''(t_2)} - \frac{Q''(0)}{2\eta''(0)^2}} \bar{w}^2 \,dy
\\ &\leq \int \paren{\frac{Q''(t_1)}{2\eta''(t_2)} |a-1| + \abs{\frac{Q''(t_1)}{2\eta''(t_2)} - \frac{Q''(0)}{2\eta''(0)^2}}} \bar{w}^2 \,dy
\\ &\leq \int \paren{\lambda \frac{Q''(t_1)}{2\eta''(t_2)} + \abs{\frac{Q''(t_1)}{2\eta''(t_2)} - \frac{Q''(0)}{2\eta''(0)^2}}} \bar{w}^2 \,dy
\\ &\leq C (\lambda + \delta) \int \bar{w}^2 \,dy
\end{aligned} \end{equation}
for $\eps_1$ and $h_1$ sufficiently small.

The $\bar{w}^2$ term is an error term: 
\begin{equation} \label{w(ubar) bound}
\frac{\lambda}{\eps} \int \frac{Q'(s)}{2 \eta''(s)} \bar{w}^2 \,dy \leq \lambda \frac{\Lambda^2}{2} \int \bar{w}^2 \,dy, 
\end{equation}
as is the $\eta'(\bar{u}|s)$ term:  by Lemma \ref{thm:taylor} \ref{taylor eta' thrd}, for $\eps_0$ and $h_0$ sufficiently small
\begin{equation} \label{eta'(ubar) bound} 
\int a \frac{Q'(s)}{\eta''(s)} \eta'(\bar{u}|s) \,dy \leq C \eps \int \bar{w}^2 \,dy.
\end{equation}
Note that $C$ here depends on $\eta'''(0)$. 

To bound the $F$ term of $B$, we utilize the formula, valid for any $f$ with $f(0)=f'(0)=0$, 
\[ f(x) = \int_0^x f''(t) (x-t) \,dt. \]
Since
\[ \frac{d^2}{dx^2} F(x;s) = Q''(x) [\eta'(x)-\eta'(s)] + Q'(x) \eta''(x), \]
we have, letting $c_Q \in [-\eps,\eps]$ be the unique point such that $Q'(c_Q) = 0$,
\begin{equation} \label{F decomposition} \begin{aligned} 
F(x;s) &= \int_s^x (x-\tau) [\eta'(\tau)-\eta'(s)]Q''(\tau) + (x-\tau) \eta''(\tau)Q'(\tau)\,d\tau
\\ &= \int_s^x (x-\tau) (\tau-s) \eta''(t_1) Q''(\tau) \,d\tau + \int_s^x \eta''(\tau) Q''(t_2) (x-\tau)(\tau-c_Q) \,d\tau
\\ &= \int_s^x \eta''(t_1) Q''(\tau) (x-\tau) (\tau-s) \,d\tau + \int_s^x \eta''(\tau) Q''(t_2) (x-\tau) (\tau-s) \,d\tau 
\\ & \hspace{4cm} + \int_s^x \eta''(\tau) Q''(t_2) (x-\tau) (s-c_Q) \,d\tau
\end{aligned} \end{equation}
for some points $t_1 \in [s,\tau]$ and $t_2 \in [c_Q,\tau]$ depending on $\tau$.  

We can estimate each of these three integrals:
\begin{equation} \label{bound the three parts} \begin{aligned}
\abs{\int_s^x \eta''(t_1) Q''(\tau) (x\!-\!\tau) (\tau\!-\!s) \,d\tau - \eta''(0)Q''(0) \frac{(x-s)^3}{6}} &\leq \sup_\tau \abs{\eta''(t_1)Q''(\tau) - \eta''(0)Q''(0)} \frac{|x-s|^3}{6}, \\
\abs{\int_s^x \eta''(\tau) Q''(t_2) (x\!-\!\tau) (\tau\!-\!s) \,d\tau - \eta''(0)Q''(0) \frac{(x-s)^3}{6}} &\leq 
\sup_\tau \abs{\eta''(\tau)Q''(t_2) - \eta''(0)Q''(0)} \frac{|x-s|^3}{6}, \\
\abs{\int_s^x \eta''(\tau) Q''(t_2) (x-\tau) (s-c_Q) \,d\tau} &\leq 2 \eps \Lambda^2 \int_x^s |x-\tau| \,d\tau = \eps \Lambda^2 (x-s)^2.
\end{aligned}\end{equation}

Therefore, if $\eps_1$ and $h_1$ are sufficiently small then from \eqref{F decomposition} and \eqref{bound the three parts} we obtain
\begin{equation} \label{F(ubar) bound}
\abs{ \frac{\lambda}{\eps} \int F(\bar{u};s) \,dy - \frac{\lambda}{\eps}\frac{Q''(0)}{3 \eta''(0)^2} \int \bar{w}^3 \,dy } \leq C \frac{\lambda}{\eps} \delta \int |\bar{w}|^3 \,dy + C \lambda \int \bar{w}^2 \,dy. 
\end{equation}

Combining \eqref{Q(ubar) bound}, \eqref{w(ubar) bound}, \eqref{eta'(ubar) bound}, and \eqref{F(ubar) bound},
\begin{equation} \label{B(ubar) bound}
B \leq \frac{\lambda}{\eps}\frac{Q''(0)}{3 \eta''(0)^2} \int \bar{w}^3 \,dy + \delta \frac{\lambda}{\eps} C \int |\bar{w}|^3 \,dy + \frac{Q''(0)}{2 \eta''(0)^2} \int \bar{w}^2 \,dy + C (\lambda + \delta + \eps) \int \bar{w}^2 \,dy. 
\end{equation}

\vskip.3cm

Lastly, we bound the quantity $D$.  Define $y_{\pm} := \eta'(\mp \eps)$.  Applying Lemma~\ref{thm:taylor} \ref{taylor 1-y2}, 
\begin{equation}\begin{aligned} \label{D(ubar) bound}
(1-h) D(\bar{u}) &\geq \frac{Q''(t_1)}{2\eta''(t_2)\eta''(t_3)} (1-h) \int [y-y_-][y_+-y] |\del_y \bar{w}|^2 \,dy
\\ &\geq \frac{Q''(0)}{2 \eta''(0)^2} (1-\delta) \int [y-y_-][y_+-y] |\del_y \bar{w}|^2 \,dy
\end{aligned}\end{equation}
so long as $\eps_1$ and $h_1$ are sufficiently small.  

\vskip.3cm

We can now bound the quantity $\bar{R}$.  By combining the bounds \eqref{Y2(ubar) bound}, \eqref{B(ubar) bound}, and \eqref{D(ubar) bound} on $Y$, $B$, and $D$ respectively,
\begin{equation} \begin{aligned} \label{R(ubar) bound}
\bar{R} &\leq \frac{-C}{\eps \delta} \paren{\frac{\lambda}{\eps} \int \bar{w}^2 \,dy + 2 \int \bar{w} \,dy}^2 + \frac{\lambda}{\eps} \frac{Q''(0)}{3 \eta''(0)^2} \int \bar{w}^3 + \frac{Q''(0)}{2 \eta''(0)^2} \int \bar{w}^2 
\\ & \hspace{1cm} - \frac{Q''(0)}{2 \eta''(0)^2} (1-\delta) \int [y-y_-][y_+-y] |\del_y \bar{w}|^2 \,dy
\\ & \hspace{2cm} + C\paren{\frac{\lambda^2}{\delta} + \lambda + \delta + \eps} \int \bar{w}^2 \,dy + C \frac{\lambda \delta}{\eps} \int |\bar{w}|^3 \,dy
\\ &= \frac{Q''(0)}{2 \eta''(0)^2} \bigg[\frac{-C}{\eps \delta}\paren{\frac{\lambda}{\eps} \int \bar{w}^2 \,dy + 2 \int \bar{w} \,dy}^2 + \frac{\lambda}{\eps} \frac{2}{3}\int \bar{w}^3 + \int \bar{w}^2 
\\ & \hspace{1cm} - (1-\delta) \int [y-y_-][y_+-y] |\del_y \bar{w}|^2 \,dy 
\\ & \hspace{2cm} + C\paren{\frac{\eps^2 \lambda^2}{\delta} + \lambda + \delta + \eps} \int \bar{w}^2 \,dy + C \frac{\lambda \delta}{\eps} \int |\bar{w}|^3 \,dy \bigg]
\end{aligned} \end{equation}

\vskip.3cm

We will now perform a change of coordinates.  Let $L := \eta'(\sne) - \eta'(\spo)$.  Consider $z \in [0,1]$ and define
\begin{align*}
y &:= \eta'(\sne) + z L, \\
dy &= L dz, \\
W(z) &:= \frac{\lambda}{\eps} \bar{w}(y) = \frac{\lambda}{\eps} \bar{w}\paren{\eta'(\sne) + z L}, \\
\del_z W(z) &= \frac{\lambda}{\varepsilon} L \del_y \bar{w}(y).  
\end{align*}
Note that $z=0$ corresponds to $y = \eta'(\sne)$ and $z=1$ to $y = \eta'(\spo)$.  

In these coordinates, 
\begin{align*}
\int \bar{w} \,dy &= \frac{\eps}{\lambda} L \int W \,dz, \\
\int \bar{w}^2 \,dy &= \frac{\eps^2}{\lambda^2} L \int W^2 \,dz, \\
\int \bar{w}^3 \,dy &= \frac{\eps^3}{\lambda^3} L \int W^3 \,dz, \\
\int [y-y_-][y_+-y] |\del_y \bar{w}|^2 \,dy &= \frac{\eps^2}{\lambda^2} L \int z(1-z) |\del_z W|^2 \,dz.
\end{align*}

In terms of $z$ and $W$, \eqref{R(ubar) bound} becomes
\begin{align*} 
\bar{R} &\leq \frac{L Q''(0)}{2 \eta''(0)^2} \frac{\eps^2}{\lambda^2}\bigg[\frac{-C_2 L}{\eps \delta}\paren{\int\!\! W^2 \,dz + 2 \!\!\int\!\! W \,dz}^2 + \frac{2}{3}\int\!\! W^3 \,dz + \int\!\! W^2 \,dz - (1-\delta)  \int\!\! z(1-z) |\del_z W|^2 \,dz
\\ & \hspace{2cm} + C_3 \paren{\frac{\eps^2 \lambda^2}{\delta} + \lambda + \delta + \eps} \int W^2 \,dz + C \delta \int |W|^3 \,dz \bigg]
\end{align*}

Fixing now $\delta$ so that $\delta < \frac{\delta_0}{3 C_3}$ and $\delta < C_2 \Lambda \delta_0$, then taking $\eps_1$ small enough that $C_3(\frac{\lambda^2}{\delta} + \delta + \eps + \lambda + \delta) \leq \delta_0$ and $\eps_1 < \delta$, and recalling $L/\eps \leq \Lambda$, we can bound
\begin{align*} 
\bar{R} &\leq C \frac{\eps^2}{\lambda^2} \bigg[ \frac{-1}{\delta_0} \paren{\int W^2 \,dz + 2 \int W \,dz}^2 + \frac{2}{3}\int W^3 \,dz + \delta_0 \int |W|^3 \,dz + (1+\delta_0)\int W^2 \,dz 
\\ & \hspace{2cm} - (1-\delta_0) \int z(1-z) |\del_z W|^2 \,dz \bigg]. 
\end{align*}

We can now apply Lemma~\ref{thm:poincare} and the proof is complete.  
\end{proof}

Now that we know $-\eps^{-2} Y^2 + B - D$ is non-negative for $u$ sufficiently close to $s$, we can bound the same quantity for $u$ large by decomposing into a part near $s$ and a part far away.  

\begin{proposition}[Decrease for large perturbations] \label{thm:functional far}
Let $\eta$ and $Q$ satisfy \eqref{convexity condition} for all $x \in \R$ and have $\eta'''$, $Q'''$ continuous at 0.  For any positive constant $\bar{C}$, there exists a constant $\eps_2 > 0$ such that the following holds:

Let $s$ be a stationary shock solution to \eqref{eq:main} with $\nu = 1$ and $s_{\pm}= \mp \eps$ with $0 < \eps < \eps_2$.  There exists a $\lambda > 0$ such that for all $u:\R \to \R$ such that $w:= \eta'(u)-\eta'(s)$ satisfies
\[ \int w^2 \leq \bar{C} \frac{\eps^3}{\lambda^2}, \]
$u$ and $a := 1-\frac{\lambda}{\eps} \eta'(s)$ satisfy
\[ R(u) := \frac{-1}{2\eps^2} Y(u)^2 + B(u) - D(u) \leq -\eps_2 D(u). \]
\end{proposition}

\begin{proof}
Let $h_1$ and $\eps_1$ be the parameters defined by Proposition~\ref{thm:functional near}, and define $\bar{u}$ for a parameter $0 < h < h_1$ such that
\[ \begin{cases}
\bar{u} = u & |\eta'(u)-\eta'(s)| \leq h, \\
\eta'(u)-\eta'(s) = h\sgn(u-s) & \textrm{else}.
\end{cases} \]
Then we can define $\bar{w} := \eta'(\bar{u}) - \eta'(s)$, $\tilde{w} := w - \bar{w}$, $\tilde{Y} := Y(u) - Y(\bar{u})$, $\tilde{B} := B(u) - B(\bar{u})$, and $\tilde{D} := D(u) - D(\bar{u})$.  For $\tilde{D}$ we have
\begin{equation} \label{far bound D}
\tilde{D} = \int a \indic{|w|>h} |\del_x (\eta'(u) - \eta'(s))|^2 dx. 
\end{equation}

\vskip.3cm
We will bound $\tilde{Y}$, $\tilde{B}$, and $\tilde{D}$ one at a time.  

To bound $\tilde{Y}$, we calculate
\begin{align*}
\eta(u|s) - \eta(\bar{u}|s) &= \int_{\bar{u}}^u \eta'(t)-\eta'(s) \,dt
\\ &\leq \Lambda \int_{\bar{u}}^u [t-s]\,dt
\\ &= \Lambda \bracket{\frac{t^2}{2} - t s}_{\bar{u}}^u
\\ &= \Lambda \bracket{\frac{(u-\bar{u})^2}{2} + (u-\bar{u})(\bar{u}-s)}
\\ &\leq C \paren{\tilde{w}^2 + h \tilde{w}}.
\end{align*}

Therefore
\begin{align*} 
\tilde{Y} &\leq C \frac{\lambda}{\eps} \paren{\int \tilde{w}^2 + h \int \tilde{w}} + \int \tilde{w}
\\ &\leq C \frac{\lambda}{\eps} \int \tilde{w}^2 \,dy + C \paren{\frac{\lambda h}{\eps} + 1} \int \tilde{w} \,dy.
\end{align*}

Since
\[ -Y(u)^2 \leq -Y(\bar{u})^2/2 + \tilde{Y}^2, \]
and $\int \tilde{w} \leq \eps^{1/2} \paren{\int \tilde{w}^2}^{1/2}$, and $\int \tilde{w}^2 \leq \bar{C} \eps^3/\lambda^2$, we can bound
\begin{equation} \label{far bound Y}
\frac{-1}{\eps^2} Y(u)^2 \leq \frac{-1}{2\eps^2} Y(\bar{u})^2 + C \eps \int \tilde{w}^2 \,dy + C \paren{\eps + \frac{\lambda^2 h^2}{\eps}} \int \tilde{w} \,dy.
\end{equation}

\vskip.3cm

For the $B$ term, we must assume without loss of generality that $2 \Lambda \eps \leq h_1$ (so that $Q'$ does not change sign between $\bar{u}$ and $u$).  Then we can calculate
\begin{equation}\begin{aligned} \label{far bound F}
F(u;s) - F(\bar{u};s) &= \int_{\bar{u}}^u Q'(t) [\eta'(t)-\eta'(s)] \,dt 
\\ &\leq \Lambda^2 \abs{\int_{\bar{u}}^u t [t-s] \,dt}
\\ &= \Lambda^2 \abs{\frac{t^3}{3} - \frac{t^2 s}{2}}_{\bar{u}}^u
\\ &= \Lambda^2 \abs{ \frac{(u-\bar{u})^3}{3} + (u-\bar{u})^2 \bar{u} - \frac{(u-\bar{u})^2 s}{2} + (u-\bar{u}) \bar{u}^2 - (u-\bar{u}) \bar{u} s }
\\ &\leq C \paren{ |\tilde{w}|^3 + h \tilde{w}^2 + h^2 |\tilde{w}| + \eps \tilde{w}^2 + \eps h |\tilde{w}| }. 
\end{aligned} \end{equation}

Similarly,
\begin{equation}\begin{aligned} \label{far bound Q}
Q(u|s) - Q(\bar{u}|s) &= \int_{\bar{u}}^u [Q'(t)-Q'(s)] \,dt
\\ &\leq \Lambda \int_{\bar{u}}^u [t-s] \,dt
\\ &= \Lambda \bracket{ \frac{(u-\bar{u})^2}{2} + (\bar{u}-s) (u-\bar{u})}
\\ &\leq C \paren{ \tilde{w}^2 + h \tilde{w} },
\end{aligned} \end{equation}
and
\begin{equation} \begin{aligned} \label{far bound eta'}
\eta'(u|s) - \eta'(\bar{u}|s) &= \int_{\bar{u}}^u \eta''(x)-\eta''(s) \,dx
\\ &\leq 2 \Lambda \abs{\int_{\bar{u}}^u \,dx}
\\ &\leq 2\Lambda^2 |\tilde{w}|,
\end{aligned} \end{equation}
and trivially
\begin{equation} \begin{aligned} \label{far bound w}
w^2 - \bar{w}^2 &= \tilde{w}^2 + 2 \bar{w} \tilde{w} 
\\ &\leq \tilde{w}^2 + 2 h |\tilde{w}|.
\end{aligned} \end{equation}

Combining \eqref{far bound F}, \eqref{far bound Q}, \eqref{far bound eta'}, and \eqref{far bound w}, we can bound $\tilde{B}$
\begin{equation}\begin{aligned} \label{far bound B}
|\tilde{B}| &\leq C \paren{\frac{\lambda}{\eps} |\tilde{w}|^3 + \bracket{\frac{\lambda}{\eps}(h + \eps) + 1 + \lambda} \tilde{w}^2 + \bracket{\frac{\lambda}{\eps} (h^2 + \eps h) + (\eps + h) + \eps + \lambda h} |\tilde{w}| }
\\ &= C\paren{\frac{\lambda}{\eps} |\tilde{w}|^3 + \bracket{ \frac{\lambda h}{\eps} + 1 + \lambda} \tilde{w}^2 + \bracket{\frac{\lambda h^2}{\eps} + \eps + h + \lambda h } |\tilde{w}| }.
\end{aligned} \end{equation}

Using \eqref{far bound D}, \eqref{far bound Y}, and \eqref{far bound B}, we can decompose the quantity $R(u)$ as
\begin{equation} \begin{aligned} \label{R decomposition}
R(u) &\leq \bracket{\frac{-1}{2\eps^2} Y(\bar{u})^2 - B(\bar{u}) - (1-h)D(\bar{u})} 
\\ &+ \bracket{\frac{1}{\eps^2} \paren{ \frac{\lambda}{\eps} \int \tilde{w}^2 \,dy  + \paren{1+\frac{\lambda h}{\eps}}\int \tilde{w} \,dy}^2 + \frac{\lambda}{\eps} \int \tilde{w}^3 \,dy + \frac{\lambda h}{\eps} \int \tilde{w}^2 \,dy - (1-h) \tilde{D}} 
\\ &+ \bracket{\frac{\lambda h^2}{\eps} \int \tilde{w} \,dy - \frac{h}{2} D(u)} - \frac{h}{2} D(u)
\\ &:= R_1 + R_2 + R_3 - \frac{h}{2} D(u).  
\end{aligned} \end{equation}

By Proposition \ref{thm:functional near}, we know $R_1 \leq 0$.  It remains to show the same for $R_2$ and $R_3$.  

\vskip.3cm

Using the fact that $\int \tilde{w}^2 \,dy \leq \bar{C} \eps^3/\lambda^2$, we can bound the quantity $R_2$
\begin{equation} \begin{aligned} \label{R2 bound} 
R_2 &\leq \frac{1}{\eps^2} \bracket{\frac{\lambda^2}{\eps^2} \paren{\int \tilde{w}^2 \,dy} \int \tilde{w}^2 \,dy + \paren{1 + \frac{h\lambda}{\eps}}^2 \paren{\int \tilde{w} \,dy}^2 } + \frac{h\lambda}{\eps} \int \tilde{w}^2 \,dy + \frac{\lambda}{\eps} \int \tilde{w}^3 \,dy - (1-h) \tilde{D}
\\ &\leq \paren{\frac{1}{\eps} + \frac{h \lambda}{\eps}} \int \tilde{w}^2 \,dy + \paren{\frac{1}{\eps^2} + \frac{h^2 \lambda^2}{\eps^4}} \paren{\int \tilde{w} \,dy}^2 + \frac{\lambda}{\eps} \int \tilde{w}^3 \,dy - (1-h) \tilde{D}. 
\end{aligned} \end{equation}

By Lemma \ref{thm:gagliardo}, we know that for any exponent $q \in (0,1)$ we have
\begin{align*}
\int \tilde{w} \,dy &\leq C_q \paren{\frac{\eps^3}{h^2\lambda^2}}^q \eps^{-1/2} \tilde{D}^{1/2}, \\
\int \tilde{w}^2 \,dy &\leq C_q \paren{\frac{\eps^3}{h^2\lambda^2}}^q \eps\n \tilde{D}, \\
\int \tilde{w}^3 \,dy &\leq C_q \paren{\int \tilde{w}^2 \,dy}^{1/2} \paren{\int \tilde{w}^4 \,dy}^{1/2} 
\leq \paren{\frac{\eps^3}{h^2\lambda^2}}^{q/2} \frac{\eps^{1/2}}{h \lambda} \tilde{D}.
\end{align*}

From these estimates with the appropriate $q$, we find that if $\eps$, $\lambda$ and $h$ are appropriately small (specifically if $\eps \leq C_h \lambda^3$ for constant $C_h$ depending on $h$) then
\[ \paren{\frac{1}{\eps} + \frac{h \lambda}{\eps}} \int \tilde{w}^2 \,dy + \paren{\frac{1}{\eps^2} + \frac{h^2 \lambda^2}{\eps^4}} \paren{\int \tilde{w} \,dy}^2 + \frac{\lambda}{\eps} \int \tilde{w}^3 \,dy \leq \frac{1}{2} \tilde{D}. \]

Plugging this estimate into \eqref{R2 bound}, and assuming without loss of generality $h < 1/2$, the quantity $R_2$ will be non-positive.  

\vskip.3cm

It remains to bound the quantity $R_3$.  

Let $f := \paren{|w| - \frac{h}{2}}_+$.  Then $\tilde{w} = \paren{f - \frac{h}{2}}_+$.  

By Lemma \ref{thm:gagliardo} with exponent $3/4$,
\[ \int f^2 \,dy \leq C \frac{\eps^{5/4}}{h^{3/2} \lambda^{3/2}} D(u). \]

By Chebyshev's inquality,
\[ \int |\tilde{w}| \,dy = \int \paren{|w|-h}_+ \,dy \leq \frac{2}{h} \int f^2 \,dy. \]

Therefore, 
\[ R_3 \leq \paren{C\frac{\lambda h^2}{\eps} \frac{2}{h} \frac{\eps^{5/4}}{h^{3/2} \lambda^{3/2}} - \frac{h}{2}} D(u) = \paren{C\frac{\eps^{1/4}}{h^{3/2}\lambda^{1/2}} - 1} \frac{h}{2} D(u). \]
So long as $\eps < \paren{C\n h^{3/2} \lambda^{1/2}}^4$, the quantity is non-positive.  

\vskip.3cm

Since $R_1$, $R_2$, and $R_3$ are all non-positive, by \eqref{R decomposition} we know $R(u) \leq -h/2 D(u) \leq -\eps_2 D(u)$.  

\end{proof}

\section{Proof of Main Theorem} \label{sec:main thm}

We will now prove Theorem~\ref{thm:main technical}.  The idea of the proof is to define the shift function $\gamma$ such that when $|Y(u)|$ is large, the $\dot{\gamma} Y$ term is negative and dominant, while when $|Y(u)|$ is small we can apply Proposition~\ref{thm:functional far}.  
\begin{proof}
Take $\eps_0$ to be the constant $\eps_2$ defined in Proposition~\ref{thm:functional far}.  

We must construct a shift function $\gamma$, so we begin by making elementary bounds on the term $B$.  Note that $u(x)-s(x)$ is guaranteed to be in $L^2$ for short time by the basic existence theorems of, for example, \cite{Se}. Moreover,
\begin{equation} \label{continuity in gamma}
\int |s(x) - s(x-\xi)|^2 \,dx \leq C (1+\sqrt(\xi)). 
\end{equation}

From the estimates of Lemma~\ref{thm:taylor}, we know that for some constant $C$,
\[ |B(u)| \leq C(\eps, \lambda, \Lambda) \paren{\int w^3 \,dy + \int w^2 \,dy + \int w \,dy}. \]
Moreover, since by H\"{o}lder's inequality $\int w^3 \,dy \leq \paren{\int w^2 \,dy}^{3/4} \paren{\int w^6 \,dy}^{1/4}$, we can further say by Lemma~\ref{thm:gagliardo}, by taking $h^2 = \frac{2 \Lambda}{\eps} \int w^2 \,dy$, that
\begin{align*} 
\int w^3 \,dy &\leq C \paren{\int w^2 \,dy}^{3/4} \paren{ \Lambda h \eps + \eps\n D^3}^{1/4}
\\ &\leq C(\eps) \paren{\int w^2 \,dy}^{7/8} + C(\eps)\paren{\int w^2 \,dy}^{3/4} D^{3/4}.
\end{align*}
It follows that
\begin{equation} \label{nonlinear gamma lipschitz}
2|B|-(1-\eps_0) D \leq C(\eps) \bracket{1 + \paren{\int w^2 \,dy}^3}. 
\end{equation}
Of course, $\int w^2 \,dy$ depends on $\gamma$.  

Define
\[ \Phi_\varepsilon(y) := 
\begin{cases} 
1 & y\leq -\varepsilon^2 \\
\frac{-y}{\varepsilon^2} & |y|\leq \varepsilon^2 \\
-1 & y\geq \varepsilon^2.
\end{cases} \]

We define $\gamma(t)$ as the solution of the nonlinear ODE:
\[ \begin{cases}
\dot{\gamma}(t) &= \Phi_\varepsilon(Y(u^\gamma)) \paren{\frac{1}{\eps^2} \paren{2|B_\gamma(u)| - (1-\eps_1) D_\gamma(u)}_+ + 1} \\
\gamma(0) &= 0
\end{cases} \]
From \eqref{continuity in gamma} and \eqref{nonlinear gamma lipschitz}, we know that
\[ \paren{2|B_\gamma(u)| - (1-\eps_1) D_\gamma(u)}_+ \leq C(\eps, \int |u(x) - s(x)|^2 \,dx) \bracket{1 + |\gamma(t)|^{3/2}}. \]
Therefore the quantity $\gamma$ exists for a short time.  

If $|Y| \geq \varepsilon^2$ then 
\begin{align*} 
\dot{\gamma} Y + B - D &\leq - 2 \paren{2|B| - (1-\eps_0)D}_+ + 1 + B - D
\\ &\leq -2|B| + (1-\eps_0)D - \eps^2 + B - D < -\eps_0 D.
\end{align*}

Alternatively, if $|Y| \leq \varepsilon^2$, then 
\[ \dot{\gamma} Y \leq - \frac{1}{\eps^2} Y^2. \]
We can therefore apply Proposition~\ref{thm:functional far} and conclude that
\[ \dot{\gamma} Y + B - D < -\eps_0 D. \]

It follows, from Proposition~\ref{thm:time derivative}, that $\int |u(x) - s(x-\gamma(t))|^2 \,dx$ is uniformly bounded so long as $\gamma$ exists.  

Now that we have a uniform bound on $\int w^2 \,dy$, the bound \eqref{nonlinear gamma lipschitz} shows that $\gamma$ exists and is Lipschitz for all time.  
\end{proof}

Lastly we prove Theorem \ref{thm:main simple}.  

\begin{proof}
The proof is by application of Theorem \ref{thm:main technical}.  

If $s$ is not of the form required by Theorem \ref{thm:main technical}, we can replace $Q$ by
\[ \tilde{Q}(x) := Q(x-a) + bx + c \]
for suitable constants $a$, $b$, and $c$ so that $s$ is stationary and centered about 0.  Recall that by the Rankine-Hugoniot condition, if $Q(\sne) = Q(\spo)$ then $s$ is stationary.  

If $\eta$ and $Q$ only satisfy the bound \eqref{convexity condition} on a compact interval $[-R,R]$ then, so long as $\norm{u}_\infty < R$, we can modify $\eta$ and $Q$ outside this region and $u$ will solve the modified \eqref{eq:main}.  

If $\nu \neq 1$, we merely consider 
\[ \tilde{u}(t,x) := u(x/\nu, t/\nu) \]
and 
\[ \tilde{s}(x) := s(x/\nu). \]
Then $\tilde{u}$ solves (and $\tilde{s}$ is a shock solution to)
\[ \frac{1}{\nu} \del_t u + \frac{1}{\nu} \del_x Q(u) = \frac{1}{\nu^2} \nu \del_{xx} \eta'(u) \]
which is equivalent to \eqref{eq:main} with $\nu = 1$.  
\end{proof}

\bibliographystyle{alpha}
\bibliography{ConservationContraction}

\end{document}